\documentclass[12pt]{article}

\def\Extended{extended}
\usepackage{hyperref}

\usepackage{graphicx, color}
\usepackage{amsmath, amssymb}
\makeatletter
\def\amsbb{\use@mathgroup \M@U \symAMSb}
\makeatother
\usepackage{bbold}
\newcommand{\bbd}[1]{{\mathbb{#1}}}
\usepackage{bbm}
\usepackage{dsfont}
\usepackage{bbold}
\usepackage{mathtools}
\usepackage[normalem]{ulem}
\usepackage{setspace}

\newif\ifextended
\ifx\Extended\Kaba\else\extendedtrue\fi
\newif\ifJapanese
\newif\iftesting
\ifextended\testingtrue\fi

\usepackage{theorem}
\theorembodyfont{\ifJapanese\rm\else\it\fi}
\newcount\minute	
\newcount\hour		
\newcount\hourMins  
\ifJapanese
\def\today%
{
  \the\year 年\,\zeroPadTwo{\the\month}月\,\zeroPadTwo{\the\day}日%
}
\fi
\def\now%
{
  \minute=\time    
  \hour=\time \divide \hour by 60 
  \hourMins=\hour \multiply\hourMins by 60
  \advance\minute by -\hourMins 
  \zeroPadTwo{\the\hour}:\zeroPadTwo{\the\minute}%
}
\def\zeroPadTwo#1%
{
  \ifnum #1<10 0\fi    
  #1
}
\title{\bf On geometrical characterizations of $\reals$-linear mappings}
\author{\bf Saka\'e Fuchino}

\date{}

\renewcommand{\baselinestretch}{1.2}
\renewcommand{\thefootnote}{(\arabic{footnote})\,}
\iftesting
\fi
\iftesting
\newcommand{\Label}[1]{\label{#1}\marginpar{{\renewcommand{\baselinestretch}{0.4}\tiny 
		  #1}}}
\else
\newcommand{\Label}[1]{\label{#1}}
\fi

\newcounter{frml}[section]
\newcounter{frmla}[section]
\def\thefrml{{\arabic{section}.\arabic{frml}}}
\def\thefrmla{{$\aleph$\arabic{section}.\arabic{frmla}}}
\def\frmlabel#1{\refstepcounter{frml}{\def\baka{#1}\ifx\baka\empty\else\label{#1}\fi}%
{\rm({\thefrml})\hfill\hfill\hfill}}
\def\frmlabela#1{\refstepcounter{frmla}{\def\baka{#1}\ifx\baka\empty\else\label{#1}\fi}%
{\rm({\thefrmla})\hfill\hfill\hfill}}
\def\xitem[#1]{\item[\frmlabel{#1}]\mbox{}%
	\iftesting\marginpar{{\renewcommand{%
				\baselinestretch}{0.6}\tiny#1}}\fi\ignorespaces}
\def\xitemq[#1]{\item[\frmlabel{#1}]\mbox{}%
	\ignorespaces}
\def\xitemd[#1]#2{\item[(\ref{#1})$#2$\hfill\hfill\hfill]}
\def\xitemA[#1]{\item[\frmlabela{#1}]\mbox{}%
	\iftesting\marginpar{{\renewcommand{%
				\baselinestretch}{0.6}\tiny#1}}\fi\ignorespaces}
\def\xitemsub[#1]#2{\item[\frmlabel{#1}$_{#2}$]\mbox{}%
	\iftesting\marginpar{{\renewcommand{%
				\baselinestretch}{0.6}\tiny#1}}\fi\ignorespaces}
\def\xxitem[#1][#2]{\item[(\ref{#1}{\makebox[1.4ex][c]{#2}})]\mbox{}%
	\iftesting\marginpar{{\renewcommand{%
				\baselinestretch}{0.6}\tiny\{#1\}\{#2\}}}\fi\ignorespaces}
\def\xitemof#1{{\rm({\ref{#1}})}}
\def\xitemAof#1{{\rm({\ref{#1}})}}
\def\xitemdof#1#2{$\mbox{\rm(\ref{#1})}#2$}

\newenvironment{xitemize}{\begin{list}{}{\parsep=0.5\smallskipamount%
			\itemindent=-0.4ex%
			\itemsep=0.5\smallskipamount\leftmargin=4em\labelwidth=3em\labelsep=0.7em}}%
							 {\end{list}}
\def\assert#1{\noindent\makebox[4.8ex][r]{\rm(\makebox[2.2ex][c]{#1})}\ \ \ignorespaces}
\def\wassert#1{\assert{#1}}
\def\wassertof#1{\makebox[4.8ex][r]{\rm(\makebox[2.2ex][c]{#1})}}%
\def\assertof#1{{\rm(#1)}}%

\def\daimaru#1{\makebox[1zw][c]{\mbox{\leavevmode\lower.08zh\hbox{%
        \rlap{\hbox to 0.76zw{\hfil\mbox{}\hfill{}\raisebox{0.03zh}{\scalebox{1.2}{○}}\hfil}}%
        \raise0.19zh\hbox to 1zw{\hfil{\hspace{0.16zw}\footnotesize#1}\hfil}}}}\,}

\newtheorem{Thm}{\ifJapanese{\bf 定理}\else {\bf Theorem}\fi}[section]

\newtheorem{Prop}[Thm]{\ifJapanese{\bf 命題}\else{\bf Proposition}\fi}
\newtheorem{Problem}[Thm]{\ifJapanese{\bf 未解決問題}\else{\bf Problem}\fi}

\newtheorem{Lemma}[Thm]{\ifJapanese{\bf 補題}\else{\bf Lemma}\fi}

\newtheorem{Cor}[Thm]{\ifJapanese{\bf 系}\else{\bf Corollary}\fi}

\newtheorem{Claim}{{\bf Claim}}[Thm]

\newtheorem{Subclaim}{{\bf Subclaim}}[Claim]

\newcommand{\prf}{\ifJapanese{\bf 証明．\ }\ignorespaces\else{\bf 
		Proof.\ \ }\ignorespaces\fi}

\newcommand{\prfofClaim}{\raisebox{-.4ex}{\Large $\vdash$\ \ }}

\newcommand{\Thmof}[1]{\ifJapanese{定理\,\ref{#1}}\else{Theorem~\ref{#1}}\fi}
\newcommand{\Problemof}[1]{\ifJapanese{未解決問題\,\ref{#1}}\else{Problem~\ref{#1}}\fi}

\newcommand{\Lemmaof}[1]{\ifJapanese{補題\,\ref{#1}}\else{Lemma~\ref{#1}}\fi}

\newcommand{\Claimof}[1]{{Claim \ref{#1}}}

\newcommand{\Thmabove}{{\ifJapanese 定理\else Theorem\fi\ \number\theThm}}

\newcommand{\Lemmaabove}{{\ifJapanese 補題\else Lemma\fi\ \number\theThm}}

\newcommand{\Subclaimabove}{{Subclaim \number\theSubclaim}}

\newsavebox{\qedbox}\sbox{\qedbox}{
{\unitlength=0.05mm \begin{picture}(40,60)
\put(0,0){\framebox(30,44)[cc]{}}
\put(30,-7){\rule{7\unitlength}{44\unitlength}}
\put(10,-7){\rule{27\unitlength}{7\unitlength}}
\end{picture}}}
\newcommand{\qed}{\mbox{}\hfill\usebox{\qedbox}}
\newcommand{\smallqed}%
{\mbox{}\smallskip\hfill\raisebox{-.4ex}{\Large $\dashv$}}
\newcommand{\qedof}[1]%
{\mbox{} \hspace*{\fill}{\usebox{\qedbox}{\tiny~(#1)}}}
\newcommand{\Qedof}[1]%
{\mbox{} \hspace*{\fill}{\usebox{\qedbox}%
{\tiny~(#1~\number\theThm)}}}
\newcommand{\QedAof}[1]%
{\mbox{} \hspace*{\fill}{\usebox{\qedbox}%
{\tiny~(#1~\number\theThmA)}}}
\newcommand{\qedofThm}{\Qedof{\ifJapanese 定理\else Theorem\fi}}

\newcommand{\qedofCor}{\Qedof{\ifJapanese 系\else Corollary\fi}}
\newcommand{\qedofProp}{\Qedof{\ifJapanese 命題\else Proposition\fi}}
\newcommand{\qedofLemma}{\Qedof{\ifJapanese 補題\else Lemma\fi}}

\newcommand{\qedskip}{\medskip}

\newcommand{\qedofClaim}%
{\mbox{}\hfill\raisebox{-.4ex}{\Large $\dashv$ }\nolinebreak%
\mbox{\tiny~(Claim~\number\theClaim)}}
\newcommand{\qedofClaimA}%
{\mbox{}\hfill\raisebox{-.4ex}{\Large $\dashv$ }\nolinebreak%
\mbox{\tiny~(Claim~A\,\number\theClaimA)}}
\newcommand{\qedofClaimAof}[1]%
{\mbox{}\hfill\raisebox{-.4ex}{\Large $\dashv$ }\nolinebreak%
\mbox{\tiny~(Claim~A\,\ref{#1})}}
\newcommand{\qedofSubclaim}%
{\mbox{}\hfill\raisebox{-.4ex}{\Large $\dashv$ }\nolinebreak%
\mbox{\tiny~(Subclaim~\number\theSubclaim)}}

\newcommand{\setof}[2]{\{#1\,:\,#2\}}
\newcommand{\ssetof}[1]{\{#1\}}

\newcommand{\subseteqand}[1]{\mathrel{\mathop{\subseteq}%
		\limits_{\scriptscriptstyle\hbox to 14pt{$\scriptscriptstyle #1$\hss}}}}

\newcommand{\mapping}[3]{#1:#2\rightarrow #3}

\newcommand{\imageof}{{}^{\,{\prime}{\prime}}}

\newcommand{\circleq}{\mathrel{{\leqslant}%
		\hspace{-0.86ex}{\lower-0.53ex\hbox{$\scriptscriptstyle\circ$}}}}

\newcommand{\restr}{\restriction}

\newcommand{\concat}{\mathop{{}^{\frown}}}
\newcommand{\droite}{\mathop{\mbox{\rule[0.64em]{1.4ex}{0.8pt}}}}
\newcommand{\interval}{\mathop{{}^{\rotatebox[origin=c]{-90}{\tt\small[}}}}
\newcommand{\reals}{{\amsbb{R}}}

\newcommand{\rationals}{{\amsbb{Q}}}

\newcommand{\ctentenc}{,{}\linebreak[0]\hspace{0.04ex}{{.}{.}{.}\hspace{0.1ex},\,}\linebreak[0]}

\newcommand{\veca}{\mathbb{a}}
\newcommand{\vecb}{\mathbb{b}}
\newcommand{\vecc}{\mathbb{c}}
\newcommand{\vecd}{\mathbb{d}}
\newcommand{\vece}{\mathbb{e}}
\newcommand{\vecp}{\mathbb{p}}
\newcommand{\vecq}{\mathbb{q}}
\newcommand{\vecx}{\mathbb{x}}

\newcommand{\zerov}{\bbd{0}}

\newcommand{\continuum}{2^{\aleph_0}}

\newcommand{\Z}{{\sf Z}}
\newcommand{\ZC}{{\sf ZC}}
\newcommand{\ZF}{{\sf ZF}}
\newcommand{\ZFC}{{\sf ZFC}}

\newcommand{\AC}{{\sf AC}}
\newcommand{\DC}{{\sf DC}}
\newcommand{\AD}{{\sf AD}}

\newcommand{\st}{such that}
\newcommand{\wrt}{with respect to}

\newcommand{\wolog}{without loss of generality}

\newcommand{\tfae}{the following are equivalent}

\newcommand{\Pkl}[2]{\ifx\bakakaba#1\bakakaba\ifx\bakakaba#2\bakakaba{\mathcal 
    P}_\kappa(\lambda)\else{\mathcal P}_\kappa(#2)\fi\else{\mathcal P}_{#1}(#2)\fi}

\newcommand{\utildeT}[1]{%
  \hbox to 0pt{\smash{$\mathop{\textstyle #1}\limits_{%
			\raisebox{0.4ex}[0pt]{$\scriptstyle\sim$}}$}\hss}%
  \relax\phantom{\mathord{{#1}_{\rule[-0.6ex]{0pt}{1pt}}}}}
\newcommand{\utildeS}[1]{%
	\hbox to 0pt{\smash{$\mathop{\scriptstyle #1}\limits_{%
				\raisebox{0.6ex}[0pt]{$\scriptscriptstyle\sim$}}$}\hss}%
	\relax\phantom{\mathord{{#1}_{\rule[-0.6ex]{0pt}{1pt}}}}}
\newcommand{\utildeSS}[1]{%
	\hbox to 0pt{$\mathop{\scriptscriptstyle #1}%
		\limits_{\scriptscriptstyle\sim}$\hss}%
		\relax\phantom{\underline{#1}}}

{\end{minipage}\end{trivlist}}

\begin{document}
\maketitle


\renewcommand{\thefootnote}{$\ast$\ }
  \footnotetext{Graduate School of System Informatics, Kobe University \\Rokko-dai 1-1, Nada, Kobe 657-8501 Japan
   \\
    \scalebox{0.95}[1]{\tt fuchino@diamond.kobe-u.ac.jp}}
\ifextended

\begin{quotation}
	\footnotesize
	\noindent
	\centerline{\normalsize\tt Contents\hspace{6em}\mbox{}}\mbox{}\\
       {\mbox{}\hspace{-1.6em}\tt\makebox[3.4ex][l]{\ref{intro}.}%
         Introduction}\ \ \dotfill\ \ {\normalsize\pageref{intro}}\\ 
       {\mbox{}\hspace{-1.6em}\tt\makebox[3.4ex][l]{\ref{lin-section-1}.}%
         Mappings which preserve proportions on each line}\ \ \dotfill\ \ {\normalsize\pageref{lin-section-1}}\\ 
       {\mbox{}\hspace{-1.6em}\tt\makebox[3.4ex][l]{\ref{lin-section-2}.}%
         Geometry of line preserving mappings}\ \ \dotfill\ \ {\normalsize\pageref{lin-section-2}}\\  
       {\mbox{}\hspace{-1.6em}\tt\makebox[3.4ex][l]{\ref{lin-section-3}.}%
         A partial answer to \Problemof{open-p}}\ \ \dotfill\ \ {\normalsize\pageref{lin-section-3}}\\  


       \noindent
       {\mbox{}\hspace{-1.6em}\tt
         References}\ \ \dotfill\ \ {\normalsize\pageref{ref}}\\ 

\end{quotation}

\fi
\renewcommand{\thefootnote}{}
\footnotetext{{\it Date:} March 9, 2020
  \qquad {\it Last update:} 
  \today\ (\now\ JST)\vspace{-1\smallskipamount}
}
\footnotetext{{\it MSC2020 Mathematical Subject Classification:}
  15A04, 51A15\vspace{-1\smallskipamount}}
\footnotetext{The research was supported by JSPS Kakenhi Grant No.20K03717.}
\ifextended
\footnotetext{\tt This is an extended version of the paper with the same title.  
\\
  All additional 
  details not to be contained in the submitted version of the paper are either typeset in 
  typewriter font (the font this paragraph is typeset) or put in separate appendices. 
  The numbering of the assertions is kept identical with the submitted version.

  The most up-to-date file of this extended version is downloadable as:
  \href{https://fuchino.ddo.jp/papers/linear-mappings-x.pdf}{\tt https://fuchino.ddo.jp/papers/linear-mappings-x.pdf}
}
\else
\footnotetext{An updated and extended version of this paper with more details and 
  proofs is downloadable as:\qquad \href{https://fuchino.ddo.jp/papers/linear-mappings-x.pdf}{\tt https://fuchino.ddo.jp/papers/linear-mappings-x.pdf}}
\fi
\renewcommand{\thefootnote}{\arabic{footnote})\,}
\section{Introduction}
\label{intro}
In this note, we consider several properties characterizing linear mappings $\mapping{f}{X}{Y}$ for
$\reals$-linear spaces $X$ and $Y$ in terms of points and lines in $\reals$-linear spaces $X$
and their images in $Y$ by the mappings.

\Thmof{P-lin-6} and its corollaries generalize the Fundamental Theorem of Affine Geometry,  
and characterize the $\reals$-linear mappings
$\mapping{f}{X}{Y}$ with $\dim(f\imageof X)>1$ as well as the corresponding affine mappings. 

The quite elementary proof of \Thmof{P-lin-6} suggests that the result might be a 
well-known fact possibly originated from the first half of the last century.
The author however could not find any appropriate reference to the result in the 
literature. 

Artstein-Avidan and Slomka \cite{artstein-avidan-slomka} lists some known variants of the 
Fundamental Theorem of affine geometry but our results do not seem to be covered by the 
assertions cited there. 

Even in the case that the 
results presented here turn out later to be a folklore, the author believes that the present 
note is written in  a sufficiently understandable manner that it can serve as a well-written expository 
article accessible for a very wide audience\footnote{Actually, the author began to write 
  the present text as an additional teaching material to his online lecture of an 
  introductory course in linear algebra when he realized that there are several statements 
  like that of \Problemof{open-p} and \Thmof{P-lin-6} whose status was not immediately 
  clear and whose treatment might slightly exceed the capacity of students in the first semester. Most of 
  the details of the following text is still kept elementary and detailed in such a way that even 
  some motivated undergraduate students should be able to follow.}. 

In the following, we shall use standard conventions and notation in set theory which might 
be slightly different from everyday mathematics.  
In particular, we distinguish sharply between elements $x\in y$ and singletons
$\ssetof{x}\subseteq y$. The image of a subset $U\subseteq X$ of the domain of a mapping
$\mapping{f}{X}{Y}$ is denoted by $f\imageof U$ with is also often denoted as $f[U]$ in the 
literature. Thus, for a closed interval $[a,b]$, we  
write $f\imageof[a,b]$ to denote the image of the interval by the function $f$ instead of 
writing $f[a,b]$ or even, more correctly, $f\left[{[a,b]}\right]$. A natural number $n$ is the set
$\ssetof{0,1\ctentenc n-1}$ and 
thus $i\in 2$ means $i=0$ or $i=1$. 

We work in \Z\ (the modern version of Zermelo's set theory without 
Axiom of Choice) if not mentioned otherwise.

\section{Mappings which preserve proportions on each line}
\Label{lin-section-1}
Let $X$ be a linear space over the scalar field $\reals$ (or $\reals$-linear space).
A subset $P\subseteq X$ is a {\it point}\/ if $P$ is a singleton. That is, if it is an 
affine subspace of $X$ of dimension $0$. A  
subset $L$ is a line if it is an affine subspace of $X$ of dimension $0$. Thus, $L\subseteq X$ is 
a line if and only if there are vectors $\veca$, $\vecd\in X$ \st\
$L=\reals\vecd+\veca=\setof{r\vecd+\veca}{r\in\reals}$. Similarly $E\subseteq X$ is a plane 
if it is an affine subspace of $X$ of dimension $2$. Two lines $L_0$, $L_1\subseteq X$ are 
parallel, if either $L_0=L_1$ or there is a plane $E\subseteq X$ with $L_0$, $L_1\subseteq E$ and
$L_0\cap L_1=\emptyset$. Note that the Parallel Postulate holds in this setting.

Recall that a mapping $\mapping{f}{X}{Y}$ for $\reals$-linear spaces $X$, $Y$ is said to be a 
{\it linear mapping}, if we have
\begin{xitemize}
\xitem[lin-1] $f(c\,\veca)=cf(\veca)$ for all $\veca\in X$ and $c\in\reals$; and
\xitem[lin-2] $f(\veca+\vecb)=f(\veca)+f(\vecb)$ for all $\veca$, $\vecb\in X$.
\end{xitemize}
To emphasize that we are talking about linear mappings between $\reals$-linear spaces, we 
shall also say $\reals$-linear mappings in contrast to the $\rationals$-linear mappings 
which are linear mappings between $\rationals$-linear spaces. 
Here, a mapping $\mapping{f}{X}{Y}$ is a $\rationals$-linear mapping, if 
\xitemdof{lin-1}{'} and \xitemof{lin-2} hold where
\begin{xitemize}
\xitemd[lin-1]{'} $f(c\,\veca)=cf(\veca)$ for all $\veca\in X$ and $c\in\rationals$.
\end{xitemize}

Note that, for any additive groups $X$, $Y$ and mapping $\mapping{f}{X}{Y}$,
if $f$ satisfies \xitemof{lin-2}, then $f(\zerov_X)=\zerov_Y$ holds, and  this fact is 
proved in most of the standard text books of linear algebra.   
\ifextended
{\tt[\,We have $f(\zerov_X)=f(\zerov_X+\zerov_X)=f(\zerov_X)+f(\zerov_X)$. Thus, by subtracting
    $f(\zerov_X)$ we obtain $\zerov_Y=f(\zerov_X)$.]}
\fi

The following Lemma is also elementary and well-known.
\begin{Lemma}
  \Label{P-lin-0} 
  Suppose that  $X$, $Y$  are $\reals$-linear spaces, and 
  $\mapping{f}{X}{Y}$ a linear mapping. Then 
  \begin{xitemize}
  \xitem[lin-3] 
    the image $f\imageof L$ of any line $L\subseteq X$ 
    is either a point or a line in $Y$; and 
  \xitem[lin-3-0] If $f\imageof L$ for a line $L\subseteq X$ is a line in $Y$ then 
    $f\restr L$ is 1-1.
  \end{xitemize}
\end{Lemma}
\prf Suppose that $L=\reals\vecd+\veca$ with $\vecd\not=\zerov_X$. Then, by linearity of $f$,
$f\imageof L=\setof{f(r\vecd+\veca)}{r\in\reals}=\setof{rf(\vecd)+f(\veca)}{r\in\reals}=\reals f(\vecd)+f(\veca)$. 
Thus, $f\imageof L=\ssetof{f(\veca)}$ if $f(\vecd)=\zerov_Y$. Otherwise, $f\imageof L$ is the 
line $\reals f(\vecd)+f(\veca)$. In the latter case, $\mapping{f\restr L}{L}{Y}$;
$r\vecd+\veca\mapsto r f(\vecd)+f(\veca)$, $r\in\reals$, is 1-1. 
\qedofLemma\qedskip

For two distinct vectors $\veca$, $\vecb\in X$, $\veca\droite\vecb$ denotes the line
$\reals(\vecb-\veca)+\veca=\setof{t(\vecb-\veca)+\veca}{t\in\reals}$.
$\veca\droite\vecb$ is the unique line $L$ 
with $\veca$, $\vecb\in L$. 

$\veca\interval\vecb$ denotes the closed interval between $\veca$ and
$\vecb$ ($\veca\interval\vecb=\setof{t(\vecb-\veca)+\veca}{0\leq t\leq 1}$).
$\veca\concat\vecb$ denotes the open interval between $\veca$ and $\vecb$:
$\veca\concat\vecb=\veca\interval\vecb\setminus\ssetof{\veca,\vecb}$. 
{\it $\vecc\in X$ divides $\veca\interval\vecb$ in ratio $r:s$} for $r$,
$s\in\reals$ with $r+s\not=0$, if $\vecc=\frac{r}{r+s}(\vecb-\veca)+\veca$.  
Note that, if $\frac{r}{r+s}\not\in[0,1]$, the point $\vecc$ dividing $\veca\interval\vecb$ 
in ratio $r:s$ is on the line $\veca\droite\vecb$ outside the 
interval $\veca\interval\vecb$ and that any point on the line $\veca\droite\vecb$ can be 
represented as a point dividing $\veca\interval\vecb$ in some ratio.

For convenience, we shall also use the notation $\veca\droite\vecb$, 
$\veca\interval\vecb$, and $\veca\concat\vecb$ for $\veca$, $\vecb$ with $\veca=\vecb$ defining 
$\veca\droite\veca=\veca\interval\veca=\ssetof{\veca}$ and $\veca\concat\vecb=\emptyset$ in 
this case.

The proof of \Lemmaof{P-lin-0} actually shows the following:
\begin{Lemma}
  \Label{P-lin-1} Suppose that  $X$, $Y$  are $\reals$-linear spaces, and 
  $\mapping{f}{X}{Y}$ a linear mapping. Then, for any $\veca$, $\vecb$, $\vecc\in X$ with
  $\veca\not=\vecb$ and $f(\veca)\not=f(\vecb)$, and for any $r$, $s\in\reals$, if 
  $\vecc$ divides $\veca\interval\vecb$ in ratio $r:s$ then $f(\vecc)$ divides
  $f(\veca)\interval f(\vecb)$ in ratio $r:s$.\\\qed
\end{Lemma}

The property of the linear mapping mentioned in \Lemmaof{P-lin-0} does not characterize 
linear mappings even when we add the condition $f(\zerov_X)=\zerov_Y$:

\begin{Lemma}
  \Label{P-lin-2} For any non zero-dimensional $\reals$-linear spaces $X$, $Y$, there are 
  non-linear 
  mappings $\mapping{f}{X}{Y}$ satisfying 
  \begin{xitemize}
  \item[\xitemof{lin-3}\ \ \ ] 
    the image $f\imageof L$ of any line $L\subseteq X$ is either a point or a line in $Y$;
  \item[\xitemof{lin-3-0}\ \ \ ] If $f\imageof L$ for a line $L\subseteq X$ is a line in $Y$ then 
    $f\restr L$ is 1-1; and 
  \xitem[lin-4] $f(\zerov_X)=\zerov_Y$. 
  \end{xitemize}
\end{Lemma}
\prf Suppose that $\vece_0\in X\setminus\ssetof{\zerov_X}$ and
$\vecd_0\in Y\setminus\ssetof{\zerov_Y}$. Let $B$ be a linear basis of $X$ with
$\vece_0\in B$. For $\vecx\in X$, let $\varphi(\vecx)\in\reals$ be the $\vece_0$-coordinate 
of $\vecx$ \wrt\ $B$. That is, let $\varphi(\vecx)\in\reals$ be \st\ there is a linear 
combination $\vecc$ of elements of $B\setminus\ssetof{\vece_0}$ \st\
$\vecx=\varphi(\vecx)\vece_0+\vecc$.

Let $\mapping{\psi}{\reals}{\reals}$ be any non-linear bijective function with 
\begin{xitemize}
\xitem[lin-5] 
  $\psi(0)=0$. 
\end{xitemize}
Let 
$\mapping{f}{X}{Y}$ be defined by $f(\vecx)=\psi(\varphi(\vecx))\vecd_0$. Clearly $f$ is 
not a linear mapping. $f(\zerov_X)=\zerov_Y$ by \xitemof{lin-5} and the definition of $f$. 

For a line $L\subseteq X$ with $L=\reals\vecb+\veca$ for some $\veca$, $\vecb\in X$, if
$\varphi(\vecb)=0$ then $f\imageof L=\ssetof{\psi(\varphi(\veca))\vecd_0}$. Otherwise, by 
bijectivity of $\psi$, we have $f\imageof L=\reals \vecd_0$ and $f\restr L$ is 1-1. 
\qedofLemma
\qedskip

In contrast, 
the property mentioned in \Lemmaof{P-lin-1} {\it does} characterize affine mappings. In particular, if 
we add the condition \xitemof{lin-4}, then we obtain a characterization of linear 
mappings. 

\begin{Lemma}
  \Label{P-lin-3} For any $\reals$-linear spaces $X$, $Y$, suppose that $\mapping{f}{X}{Y}$ 
  is \st\ $f$ satisfies 
  \begin{xitemize}
    \item[\xitemof{lin-4}\ \ \ ] $f(\zerov_X)=\zerov_Y$; 
  \xitem[lin-6] for any $\veca$, $\vecb\in X$ with $\veca\not=\vecb$ and
    $f(\veca)\not=f(\vecb)$, $\vecc\in X$, and $r$, $s\in\reals$, if $\vecc$ divides $\veca\interval\vecb$ in 
    ratio $r:s$, then $f(\vecc)$ divides $f(\veca)\interval f(\vecb)$ in $r:s$.
  \end{xitemize}
  Then $f$ is a linear mapping. 
\end{Lemma}
\prf Suppose that $\mapping{f}{X}{Y}$ satisfies \xitemof{lin-4} and \xitemof{lin-6}. If
$f\imageof X=\ssetof{\zerov_Y}$, then $f$ is a linear mapping. Thus, let us assume that there is
$\vecb\in X$ \st\ $f(\vecb)\not=\zerov_Y$. $\vecb\not=\zerov_X$ by \xitemof{lin-4}.

For
$\vecc\in\reals\vecb$, if $\vecc=c\,\vecb$ for some $c\in\reals$,  then $\vecc$ divides $\zerov_X\interval\vecb$ in 
ratio $c:1-c$,  
since $\vecc=c\,\vecb=\frac{c}{c+(1-c)}(\vecb-\zerov_X)+\zerov_X$. By \xitemof{lin-6}, it 
follows that $f(c\,\vecb)$ divides $f(\zerov_X)\interval f(\vecb)$ in ratio $c:1-c$. That is, 
$f(c\,\vecb)
=\underbrace{\textstyle\frac{c}{c+(1-c)}}_{=1}(f(\vecb)-\underbrace{f(\zerov_X)}_{=\zerov_Y})
+\underbrace{f(\zerov_X)}_{=\zerov_Y}=c f(\vecb)$. Together with \xitemof{lin-4}, this 
implies that $f$ satisfies \xitemof{lin-1}.

If $\veca=\vecb$, then
$f(\veca+\vecb)=f(2\veca)\underbrace{=}_{\mbox{\footnotesize by \xitemof{lin-1}}}2f(\veca)
=f(\veca)+f(\veca)=f(\veca)+f(\vecb)$.

If $\veca\not=\vecb$, then $\veca+\vecb$ divides $2\veca\interval 2\vecb$ in ratio $1:1$. 
Thus $f(\veca+\vecb)$ divides
$f(2\veca)\interval f(2\vecb)\underbrace{=}_{\mbox{\footnotesize by \xitemof{lin-1}}}2f(\veca)\interval 2f(\vecb)$ 
in ratio 1:1 by \xitemof{lin-6}. This  
means that $f(\veca+\vecb)=\frac{1}{1+1}(2f(\vecb)-2f(\veca))+2f(\veca)=f(\veca)+f(\vecb)$. 
This shows that $f$ also satisfies \xitemof{lin-2}.
\qedofLemma
\qedskip

\Lemmaof{P-lin-2} and \Lemmaof{P-lin-3} suggest the following question.

\begin{Problem}
  \Label{open-p}
  Suppose that $X$, $Y$ are $\reals$-linear spaces and $\mapping{f}{X}{Y}$ is \st
  \begin{xitemize}
  \item[\xitemof{lin-3}\ \ \ ]  the image $f\imageof L$ of any line $L\subseteq X$ 
    is either a point or a line in $Y$; 
  \item[\xitemof{lin-4}\ \ \ ] $f(\zerov_X)=\zerov_Y$; and 
  \xitem[lin-8] there are $\veca_0$, $\veca_1\in X$ \st\ $f(\veca_0)$ and $f(\veca_1)$ are linearly independent. 
  \end{xitemize}
  Does it follow that $f$ is an $\reals$-linear mapping?
\end{Problem}

In the next sections, we shall prove a characterization of linear mappings
$\mapping{f}{X}{Y}$ satisfying \xitemof{lin-8} which is slightly 
stronger than the conditions in \Problemof{open-p} (see \Thmof{P-lin-6}). 

A mapping $\mapping{f}{X}{Y}$ on linear spaces $X$, $Y$ is said to be an {\it additive} function if 
$f$ satisfies \xitemof{lin-2}. Note that, for $\reals$-linear spaces $X$, $Y$,
$\mapping{f}{X}{Y}$ is additive if and only if it is $\rationals$-linear.

Whether all additive functions $\mapping{f}{\reals}{\reals}$ are 
linear mappings is a question whose answer depends on the axioms of set-theory. In the 
Zermelo's axiom system \ZC\ of set theory with full Axiom of Choice (\AC), we can 
use a Hamel basis of $\reals$ over $\rationals$, to 
construct $2^{\continuum}$ many $\rationals$-linear (and hence additive) functions from
$\reals$ to $\reals$. Since there are only $2^{\aleph_0}$ many $\reals$-linear functions 
from $\reals$ to $\reals$, there are $2^{\continuum}$ many additive functions from $\reals$ 
to $\reals$ which are not $\reals$-linear. Since a continuous function
$\mapping{f}{\reals}{\reals}$ is decided by the information on $f\restr\rationals$, the 
next Lemma follows:
\begin{Lemma}
  \Label{P-lin-4} Suppose that $\mapping{f}{\reals}{\reals}$ is an additive function. 
  Then the following are equivalent\,: \assertof{a} $f$ is $\reals$-linear. 
  \assertof{b} $f$ is continuous. \assertof{c} $f$ is monotonous.\qed
\end{Lemma}

The equivalence of \assertof{a} $\Leftrightarrow$ \assertof{b} and \assertof{c} above 
follows from the fact that \assertof{1} if an additive function is discontinuous at a 
point, then it is everywhere discontinuous; and \assertof{2} a monotone function can be 
discontinuous at most at countably many points. 

More generally, for mappings $\mapping{f}{\reals^m}{\reals^n}$, the equivalence of 
\assertof{a} and \assertof{b} in the previous Lemma still holds. Further, 
by theorems of Steinhaus and Kuratowski, the continuity of $f$ in \assertof{b} of 
\Lemmaabove\ in this generalized setting can be replaced either 
by being a Baire function or by being a measurable function (Kuczma \cite{kuczma}. 
\cite{fuchino} contains a slightly simplified the proof of these results). 
Thus we obtain 
\begin{Thm}
  \Label{P-lin-4-a-0}
  For any $\reals$-linear spaces $X$, $Y$ \st\ $Y$ has a linear base \footnote{Of course,  
    under \AC, this extra assumption of the existence of a linear base is superfluous.}, and for any aditive
  $\mapping{f}{X}{Y}$,  
  \tfae\/:
  \assertof{a} $f$ is a linear mapping; \assertof{b} $f$ is continuous; \assertof{c} $f$ is 
  Lebesgue measurable; \assertof{d} $f$ is a Baire function. 
\end{Thm}

Solovay \cite{solovay} constructed, arguing in the theory \ZFC\ $+$ ``there is an 
inaccessible cardinal'', a model of \ZF\ together with the Axiom of Dependent Choice (\DC, a 
weakening of the full \AC\ which covers most of the usages of \AC\ in everyday mathematics) which 
also satisfies the statement that ``all subsets of $\reals$  
have Baire property and are Lebesgue measurable''. Note that, in this model, all functions 
$\mapping{f}{\reals^m}{\reals^n}$ are Baire functions and also measurable, and thus all 
additive functions from an $\reals$-linear space $X$ to an $\real$-linear space $Y$ with a 
linear base are linear.

Shelah \cite{shelah} proved that the statement ``all subsets of $\reals$ are Lebesgue 
measurable'' is equiconsistent with \ZF\ $+$ ``there is an 
inaccessible cardinal''. Thus, the inaccessible cardinal in Solovay's 
model is unavoidable for the statement about Lebesgue measurability. However, Shelah 
\cite{shelah} also shows, that \ZF\ $+$ \DC\ $+$ ``all  
subsets of $\reals$ have Baire property'', hence also the theory \ZF\ $+$ \DC\ $+$ ``all 
additive function on $\reals$ is linear'' is equiconsistent with \ZF.
This is one of the most grave signs of asymmetry lying between measure and category 
in spite of the strong similarity between them as is seen in Oxtoby \cite{oxtoby}.

Note that Lebesgue measurability and Baire property of all functions
$\reals\rightarrow\reals$ are theorems under \ZF\ $+$ the Axiom of Determinacy (\AD). This 
is a result by Jan Mycielski and Stanis\l aw Swierczkowski for Baire property and by Banach, 
Mazur for Lebesgue measurability. Thus, the assertion ``all 
additive function from an $\reals$-linear space to an $\reals$-linear space with a linear 
base is a linear mapping'' is a theorem under this axiom system.

The consistency strength of the Axiom of Determinacy, however, is much higher than that of the statement 
that all subsets of $\reals$ are Lebesgue measurable: Woodin's famous theorem states that 
the Axiom of Determinacy (over \ZF) is equiconsistent with infinitely many Woodin 
cardinals which is much stronger than, e.g. class many measurable cardinals. 

Although, in the following, we never rely on these deep results in set-theory, they remind the possibility 
that a set-theoretic subtlety of this type can lurk anywhere in a discussion like the 
following, and warns us that careful examination is necessary. 

The next lemma is an easy application of \Lemmaof{P-lin-4}:
\begin{Lemma}{\rm (Kuczma \cite{kuczma}, Theorem 14.4.1)} 
  \Label{P-lin-4-a} suppose that $\mapping{f}{\reals}{\reals}$ is an additive function and 
  also multiplicative, i.e. we have that 
  \begin{xitemize}
  \xitem[lin-8-0] $f(rs)=f(r)f(s)$ for all $r$, $s\in\reals$.
  \end{xitemize}
  Then $f$ is either the constant zero function (i.e.\ $f(r)=0$ for all $r\in\reals$) or the identity 
  function (i.e.\ $f(r)=r$ for all $r\in\reals$).
\end{Lemma}
\prf
For any $x\in\reals$,
\begin{xitemize}
\xitem[lin-8-1]  if $x\geq 0$, then $f(x)=f((\sqrt{x})^2)=(f(\sqrt{x}))^2\geq 0$
\end{xitemize}
by the multiplicativity \xitemof{lin-8-0}. 
It follows that,  
for any $x$, $y\in\reals$ with $x\leq y$, 
\begin{xitemize}
\xitem[lin-8-2]
 $f(y)=f(x+(y-x))=f(x)+f(y-x)\geq f(x)$
\end{xitemize}
by additivity of $f$. 
Thus $f$ is a monotone function. By \Lemmaof{P-lin-4}, and since $f$ 
is an additive function, there is $c\in\reals$ \st\ 
$f(x)=cx$ holds for all $x\in\reals$. 

By \xitemof{lin-8-0}, it follows that 
\begin{xitemize}
\xitem[lin-8-3]
  $c=f(1)=f(1\cdot 1)=f(1)f(1)=c^2$.
\end{xitemize}

Since $f(1)\geq0$ by \xitemof{lin-8-1},
it follows that $c=0$ or $c=1$. If $c=0$, $f$ is the constant function 
$f(x)=0$ for all $x\in\reals$. If $c=1$,  $f=id_\reals$. 
\qedofLemma


\begin{Thm}
  \Label{P-lin-4-0}
  Suppose that $X$ and $Y$ are $\reals$-linear spaces and $\mapping{f}{X}{Y}$ is an 
  additive function (i.e.\ it satisfies \xitemof{lin-2}). If there is a function 
  $\mapping{\varphi}{\reals}{\reals}$ with 
  \begin{xitemize}
  \xitem[lin-9-0] $f(r\veca)=\varphi(r)f(\veca)$ for all $\veca\in X$ and $r\in\reals$, 
  \end{xitemize}
  then $f$ is an $\reals$-linear mapping.\smallskip
\end{Thm}
\prf If $f\imageof X=\ssetof{\zerov_Y}$, then $f$ is a linear mapping. Thus, we may assume that
$f\imageof X\not=\ssetof{\zerov_Y}$. 

Then we have $\varphi(1)=1$.
Hence, 
by \Lemmaof{P-lin-4-a}, it is enough to show that $\varphi$ is additive and multiplicative.

To show that $\varphi$ is additive, suppose that $r$, $s\in\reals$. Let 
$\veca\in X$ be \st\ $f(\veca)\not=\zerov_Y$. By additivity of $f$, we have 
$\varphi(r+s)f(\veca)=f((r+s)\veca)=f(r\veca+s\veca)=f(r\veca)+f(s\veca)=\varphi(r)f(\veca)+\varphi(s)f(\veca)
=(\varphi(r)+\varphi(s))f(\veca)$. It follows 
that $\varphi(r+s)=\varphi(r)+\varphi(s)$.

Multiplicativity of $\varphi$ can be shown similarly: Suppose $r$, $s\in\reals$, and let
$\veca\in X$ be \st\ $f(\veca)\not=\zerov_Y$. Then we have
$\varphi(rs)f(\veca)=f(rs\veca)=\varphi(r)f(s\veca)=\varphi(r)\varphi(s)f(\veca)$. It 
follows that $\varphi(rs)=\varphi(r)\varphi(s)$. 
\qedofThm

\section{Geometry of line preserving mappings}
\Label{lin-section-2}

\begin{Lemma}
  \Label{P-lin-5} For $\reals$-linear spaces $X$, $Y$ let $\mapping{f}{X}{Y}$ be a mapping 
  \st
  \begin{xitemize}
  \item[\xitemof{lin-3}\ \ \ ]   the image $f\imageof L$ of any line $L\subseteq X$ 
    is either a point or a line in $Y$; and
  \xitem[lin-7] for any line $L\subseteq X$, if $f\imageof L$ is also a line in $Y$, then
    $f\restr L$ is 1-1. 
  \end{xitemize}
  Then, \wassertof{1} $f$ maps any plane
  $E\subseteq X$ to either a plane or a line or a point of $Y$. If $f\imageof E$ is a plane 
  then $f\restr E$ is 1-1.\smallskip

  \wassert{2} Suppose that $L_0$, $L_1\subseteq X$ are lines \st\ $f\imageof L_0$ and
  $f\imageof L_1$ are also lines in $Y$. If $L_0$ and $L_1$ are parallel to each other then
  $f\imageof L_0$ and $f\imageof L_1$ are also parallel to each other. 
\end{Lemma}
\prf \assertof{1}: Suppose that $E\subseteq X$ is a plane.
Let $\veca$, $\vecb$, $\vecc\in E$ be \st\ $\vecb-\veca$ and $\vecc-\veca$ are independent. 
\smallskip

{\bf Case 1.} $f(\veca)=f(\vecb)$ and $f(\veca)=f(\vecc)$. By \xitemof{lin-7} and \xitemof{lin-3}, we have 
$f\imageof(\veca\droite\vecb)=f\imageof(\veca\droite\vecc)=\ssetof{f(\veca)}$. Suppose
$\vecp\in E\setminus (\veca\droite\vecb\cup\veca\droite\vecc)$. Then there is a line
$L\subseteq X$ going through $\vecp$ and crossing $\veca\droite\vecb$ and
$\veca\droite\vecc$ at two different points. Since the value of $f$ at both of these two points 
is $f(\veca)$, $f\restr L$ is not 1-1. By \xitemof{lin-7} and \xitemof{lin-3} it follows 
that $f(\vecp)=f(\veca)$. Thus we have $f\imageof E=\ssetof{f(\veca)}$. \smallskip

{\bf Case 2.} $f(\veca)=f(\vecb)$ and $f(\veca)\not=f(\vecc)$. Then
$f\imageof(\veca\droite\vecb)=\ssetof{f(\veca)}$, 
$f\imageof(\veca\droite\vecc)=f(\veca)\droite f(\vecc)$, and $f\imageof(\veca\droite\vecc)$ is a line in $Y$. 
For
$\vecp\in E\setminus (\veca\droite\vecb\cup\veca\droite\vecc)$, there is a line
$L\subseteq X$ going through $\vecp$ and crossing $\veca\droite\vecb$ and
$\veca\droite\vecc$ at two different points, say at $\ssetof{\vecp_0}$ and $\ssetof{\vecp_1}$ respectively. If
$f(\vecp_0)=f(\vecp_1)$, then we have $f(\vecp)=f(\vecp_0)=f(\veca)$. Otherwise, since
$f(\vecp_0)$ and $f(\vecp_1)$ belong to $f(\veca)\droite f(\vecc)$ we have
$f(\vecp)\in f(\vecp_0)\droite f(\vecp_1)=f(\veca)\droite f(\vecc)$. This shows that
$f\imageof E=f(\veca)\droite f(\vecc)$.
\smallskip

{\bf Case 3.} $f(\veca)\not=f(\vecb)$ and $f(\veca)=f(\vecc)$. This case can be treated 
similarly to the Case 2.\ to conclude that $f\imageof E=f(\veca)\droite f(\vecb)$.
\smallskip

{\bf Case 4.} $f(\veca)\not=f(\vecb)$ and $f(\veca)\not=f(\vecc)$. Then we have 
$f\imageof(\veca\droite\vecc)=f(\veca)\droite f(\vecb)$, 
$f\imageof(\veca\droite\vecc)=f(\veca)\droite f(\vecc)$, and both $f\imageof(\veca\droite\vecb)$ and
$f\imageof(\veca\droite\vecc)$ are lines in $Y$. If
$f(\veca)\droite f(\vecb)=f(\veca)\droite f(\vecc)$ then we can argue similarly to Case 2.\ 
and show $f\imageof E=f(\veca)\droite f(\vecb)$.

Suppose now $f(\veca)\droite f(\vecb)\not=f(\veca)\droite f(\vecc)$. Then we have
$f(\veca)\droite f(\vecb)\ \cap\ f(\veca)\droite f(\vecc)\\=\ssetof{f(\veca)}$. Let $E^*$ be the 
plane with $f(\veca)\droite f(\vecb)$, $f(\veca)\droite f(\vecc)\subseteq E^*$. 

If
$\vecp\in E\setminus (\veca\droite\vecb\cup\veca\droite\vecc)$ then we can take a line
$L\subseteq X$ going through $\vecp$ and crossing $\veca\droite\vecb$ and
$\veca\droite\vecc$ at two different points, say at $\ssetof{\vecp_0}$ 
and $\ssetof{\vecp_1}$ respectively. Since $f(\vecp_0)$,
$f(\vecp_1)\in f(\veca)\droite f(\vecb)\ \cup\ f(\veca)\droite f(\vecc)$, we have
$f(\vecp)\in f(\vecp_0)\droite f(\vecp_1)\subseteq E^*$. Conversely, suppose that
$\vecq\in E^*$. Then there is a line $L^*$ in $Y$ going through $\vecq$ and crossing
$f(\veca)\droite f(\vecb)$ and $f(\veca)\droite f(\vecc)$ at two different points say
$\ssetof{\vecq_0}$ and $\ssetof{\vecq_1}$. Let $\vecp_0\in \veca\droite\vecb$,
$\vecp_1\in\veca\droite\vecc$ be the unique vectors \st\ $f(\vecp_0)=\vecq_0$ and
$f(\vecp_1)=\vecq_1$. Since $f\imageof(\vecp_0\droite\vecp_1)=\vecq_0\droite\vecq_1$, there 
is $\vecp\in\vecp_0\droite\vecp_1\subseteq E$ \st\ $f(\vecp)=\vecq$. 

To see that $f\restr E$ is 1-1, suppose that $\vecp_0$, $\vecp_1\in E$ are distinct 
vectors. If $f(\vecp_0)=f(\vecp_1)$ then we can produce the constellation of Case 1.\ or Case 
2.\ with $\vecp_0$ and $\vecp_1$ together with some third point $\in E$. This is a 
contradiction since we already know that $f\imageof E$ is a plane in $Y$. 
\smallskip

\assertof{2}: If $f\imageof L_0=f\imageof L_1$ then $f\imageof L_0$ and $f\imageof L_1$ are 
parallel.

Assume otherwise. Let $E$ be the plane in $X$ with $L_0$, $L_1\subseteq E$. Since $f\imageof L_0$,
$f\imageof L_1\subseteq f\imageof E$, and since $f\imageof L_0$ and $f\imageof L_1$ are two 
distinct lines in $f\imageof E$, $f\imageof E$ is a plane in $Y$ 
and $f\restr E$ is 1-1, by \assertof{1}. Since
$L_0\cap L_1=\emptyset$. We also have $f\imageof L_0\cap f\imageof L_1=\emptyset$. Thus, 
also in this case, $f\imageof L_0$ and $f\imageof L_1$ are parallel to each other.
\qedofLemma

\begin{Lemma}
  \Label{P-lin-5-0} Suppose that $X$ and $Y$ are $\reals$-linear spaces, and
  $\mapping{f}{X}{Y}$ a mapping satisfying \xitemof{lin-4}, \xitemof{lin-3} and 
  \xitemof{lin-7}. If $\veca$, $\vecb\in X$ are \st\ $\veca\not=\zerov_X$,
  $\vecb\not=\zerov_X$, $f(\veca)=\zerov_Y$, and $f(\vecb)\not=\zerov_Y$, then we have
  $f(\veca+\vecb)=f(\vecb)=f(\veca)+f(\vecb)$. 
\end{Lemma}
\prf $\veca$ and $\vecb$ are linearly independent: $f\restr\reals\vecb$ is 1-1 by 
\xitemof{lin-3} and \xitemof{lin-7}. Thus, if $\veca\in\reals\vecb$, we would have $f(\veca)\not=\zerov_Y$ 
by \xitemof{lin-4}.

Let $E\subseteq X$ be the plane with $\zerov_X$, $\veca$, $\vecb\in E$. 
By \Lemmaof{P-lin-5},\,\assertof{1}, we have $f\imageof E=f(\vecb)\droite\zerov_Y=\reals f(\vecb)$.
Note also that $f\imageof\reals\veca=\ssetof{\zerov_Y}$. 

Suppose, toward a contradiction, that 
\begin{xitemize}
\xitem[lin-9-1] $f(\veca+\vecb)\not=f(\vecb)$.
\end{xitemize}
Then $f\restr(\veca+\vecb)\droite\vecb$ is not a constant function. Thus, we have
$f\imageof(\veca+\vecb)\droite\vecb=\reals f(\vecb)$. 
Let $\vecc\in(\veca+\vecb)\droite\vecb$ be \st\ $f(\vecc)=\zerov_Y$. 
Then we have $f\imageof\vecc\droite\zerov_X=\ssetof{\zerov_Y}$.

For each $\vecd\in(\veca+\vecb)\droite\vecb$, let $L\subseteq X$ be a line going through 
$\vecd$ which crosses $\vecc\droite\zerov_X$ and $\veca\droite\zerov_X$ at different points, say at points
$\vecp_0$ and $\vecp_1$ respectively. Then we have $f(\vecp_0)=f(\vecp_1)=\zerov_Y$. Thus
$f\imageof L=\ssetof{\zerov_Y}$ and $f(\vecd)=\zerov_Y$. In particular, we have
$f(\veca+\vecb)=f(\vecb)=\zerov_Y$. This is a contradiction to the assumption \xitemof{lin-9-1}.\qedofLemma

\section{A partial answer to \Problemof{open-p}}
\Label{lin-section-3}

The following theorem gives a positive partial answer to \Problemof{open-p}

\begin{Thm}
  \Label{P-lin-6} Suppose that $X$, $Y$ are $\reals$-linear spaces and 
  $\mapping{f}{X}{Y}$ satisfies
  \begin{xitemize}
  \item[\xitemof{lin-4}\ \ \ ] $f(\zerov_X)=\zerov_Y$; and, 
  \item[\xitemof{lin-8}\ \ \ ] there are $\veca_0$, $\veca_1\in X$ \st\ $f(\veca_0)$ and
    $f(\veca_1)$ are linearly independent.       
  \end{xitemize}

  Then $f$ is a linear mapping if and only if $f$ satisfies
  \begin{xitemize}
  \item[\xitemof{lin-3}\ \ \ ] the image $f\imageof L$ of any line $L\subseteq X$ by $f$ 
    is either a point or a line in $Y$; and,  
  \item[\xitemof{lin-7}\ \ \ ] for any line $L\subseteq X$, if $f\imageof L$ is also a line in $Y$, then
    $f\restr L$ is 1-1.
  \end{xitemize}



\end{Thm}
\prf Linear mappings satisfy \xitemof{lin-3} and \xitemof{lin-7} (\Lemmaof{P-lin-0}).

Suppose that $\mapping{f}{X}{Y}$ satisfies \xitemof{lin-4}, \xitemof{lin-8}, 
\xitemof{lin-3} and \xitemof{lin-7}. By \Thmof{P-lin-4-0}, it is enough to show that there 
is $\mapping{\varphi}{\reals}{\reals}$ with \xitemof{lin-9-0} for this $f$,  and $f$ is additive. 

To prove that there is a mapping $\mapping{\varphi}{\reals}{\reals}$ with 
\xitemof{lin-9-0}, suppose $\veca\in X$. Note that
$\veca\droite\zerov_X=\reals\veca$. If $f(\veca)=\zerov_Y$, then 
$f\imageof(\veca\droite\zerov_X)=\ssetof{\zerov_Y}$ by
\xitemof{lin-3}, \xitemof{lin-4} and  
\xitemof{lin-7}. Thus we have
\begin{xitemize}
\xitem[lin-10] $f(r\veca)=\zerov_Y=r\,f(\veca)$, for any $r\in\reals$ (under $f(\veca)=\zerov_Y$). 
\end{xitemize}

Suppose now that $f(\veca)\not=\zerov_Y$. Then
$f\imageof\reals\veca=f(\veca)\droite\zerov_Y$
($=\reals f(\veca)$) and 
$f\restr\reals\veca$ is 1-1 by \xitemof{lin-3}, \xitemof{lin-4}, 
\xitemof{lin-7}. For each $\veca\in X$ with $f(\veca)\not=\zerov_Y$ and $r\in\reals$, let 
$\varphi_0(\veca,r)$ be $s\in\reals$ \st\ $f(r\veca)=\varphi_0(\veca,r)f(\veca)$. 
\begin{Claim}
  \Label{Cl-lin-0} $\varphi_0(\veca,r)$ does not depend on $\veca\in X$ with
  $f(\veca)\not=\zerov_Y$. 
\end{Claim}
\prfofClaim
we first prove the following:
\begin{Subclaim}
For any $\veca$, $\vecb\in X$ \st\ $f(\veca)$ and
$f(\vecb)$ are independent in $Y$,  
we have $\varphi_0(\veca,r)=\varphi_0(\vecb,r)$ for all $r\in\reals$.  
\end{Subclaim}
\prfofClaim
For $\veca\in X$ with $f(\veca)\not=\zerov_Y$, $f(0\veca)=\zerov_Y$ by \xitemof{lin-4}. 
Thus $\varphi_0(\veca,0)=0$ for all such $\veca\in X$. 

For $r\not=0$, 
let us consider the four lines
$L_0=\veca\droite\zerov_X$, $L_1=\vecb\droite\zerov_X$, $L_2=\veca\droite\vecb$, and 
$L_3=r\veca\droite r\vecb$. We have $L_0\cap L_1=\ssetof{\zerov_X}$,
$L_0\cap L_2=\ssetof{\veca}$, $L_1\cap L_2=\ssetof{\vecb}$, $L_0\cap L_3=\ssetof{r\veca}$,
$L_1\cap L_3=\ssetof{r\vecb}$, and, $L_2$ and $L_3$ are parallel to each other.  By \Lemmaof{P-lin-5}, it follows that 
$f\imageof L_0=f(\veca)\droite\zerov_Y$, $f\imageof L_1=f(\vecb)\droite\zerov_Y$,
$f\imageof L_2=f(\veca)\droite f(\vecb)$, and 
$f\imageof L_3=f(r\veca)\droite f(r\vecb)=\varphi_0(\veca,r)f(\veca)\droite\varphi_0(\vecb,r)f(\vecb)$.
$f\imageof L_0\cap f\imageof L_1=\ssetof{\zerov_Y}$,
$f\imageof L_0\cap f\imageof L_2=\ssetof{f(\veca)}$,
$f\imageof L_1\cap f\imageof L_2=\ssetof{f(\vecb)}$,
$f\imageof L_0\cap f\imageof L_3=\ssetof{f(r\veca)}=\ssetof{\varphi_0(\veca,r)f(\veca)}$, 
$f\imageof L_1\cap f\imageof L_3=\ssetof{f(r\vecb)}=\ssetof{\varphi_0(\vecb,r)f(\vecb)}$, 
and, $f\imageof L_2$ and $f\imageof L_3$ are parallel to each other.  The last condition is only possible 
when $\varphi_0(\veca,r)=\varphi_0(\vecb,r)$. 
\qedofSubclaim\qedskip

Let $\veca_0$, $\veca_1\in X$ be as in \xitemof{lin-8}. By \Subclaimabove, we have
$\varphi_0(\veca_0,r)=\varphi_0(\veca_1,r)$ for any $r\in\reals$. Suppose that $\veca\in X$ 
is \st\ $f(\veca)\not=\zerov_Y$. Then there is $i\in 2$, \st\ $f(\veca)$ and $f(\veca_i)$ 
are linearly independent. By \Subclaimabove, we 
have $\varphi_0(\veca,r)=\varphi_0(\veca_i,r)$. Thus, for any $\veca\in X$ with
$f(\veca)\not=\zerov_Y$ and $r\in\reals$, $\varphi_0(\veca,r)=\varphi_0(\veca_0,r)$. 
\qedofClaim\qedskip

Let $\mapping{\varphi}{\reals}{\reals}$ be defined by $\varphi(r)=\varphi_0(\veca_0,r)$. By 
\Claimof{Cl-lin-0} and the argument above it, this $\varphi$ satisfies \xitemof{lin-9-0}. \smallskip

To prove the additivity of $f$, suppose $\veca$, $\vecb\in X$. 
\smallskip

{\bf Case 1.} $f(\veca)$ and $f(\vecb)$ are independent in $Y$. In this case, $\veca$ and 
$\vecb$ are independent and, letting $E$ be the plane in $X$ with $\veca$, $\vecb$,
$\zerov_X\in E$, we have that $f\restr E$ is a 1-1 mapping and $f\imageof E$ is the plain in 
$Y$ with $f(\veca)$, $f(\vecb)$, $\zerov_Y$ ($=f(\zerov_X)$) $\in f\imageof E$ 
by \Lemmaof{P-lin-5},\,\assertof{1}.

Let $L_0=\veca\droite\zerov_X$, $L_1=\vecb\droite\zerov_X$,
$L_2=\veca\droite(\veca+\vecb)$, $L_3=\vecb\droite(\veca+\vecb)$. Then we have
$L_0\cap L_1=\ssetof{\zerov_X}$, $L_0\cap L_2=\ssetof{\veca}$,
$L_1\cap L_3=\ssetof{\vecb}$, and $L_2\cap L_3=\ssetof{\veca+\vecb}$. We also have that 
$L_0$ and $L_3$ are parallel to each other, and,  $L_1$ and $L_2$ are parallel to each 
other. 

By \xitemof{lin-3} and \xitemof{lin-4}, and since $f\restr E$ is 1-1, it follows 
that $f\imageof L_0=f(\veca)\droite\zerov_Y$, $f\imageof L_1=f(\vecb)\droite\zerov_Y$, 
$f\imageof L_2=f(\veca)\droite f(\veca+\vecb)$,
$f\imageof L_3=f(\vecb)\droite f(\veca+\vecb)$. We also have 
$f\imageof L_0\cap f\imageof L_1=\ssetof{\zerov_Y}$, $f\imageof L_0\cap \imageof L_2=\ssetof{f(\veca)}$,
$f\imageof L_1\cap f\imageof L_3=\ssetof{f(\vecb)}$, and $f\imageof L_2\cap f\imageof L_3=\ssetof{f(\veca+\vecb)}$.
$f\imageof L_0$ and $f\imageof L_3$ are parallel to each other, and,  $f\imageof L_1$ and
$f\imageof L_2$ are parallel to each 
other by \Lemmaof{P-lin-5},\,\assertof{2}. This implies that
$f(\veca+\vecb)=f(\veca)+f(\vecb)$.
\smallskip

\noindent
\mbox{}\hfill
\includegraphics[width=0.8\textwidth]{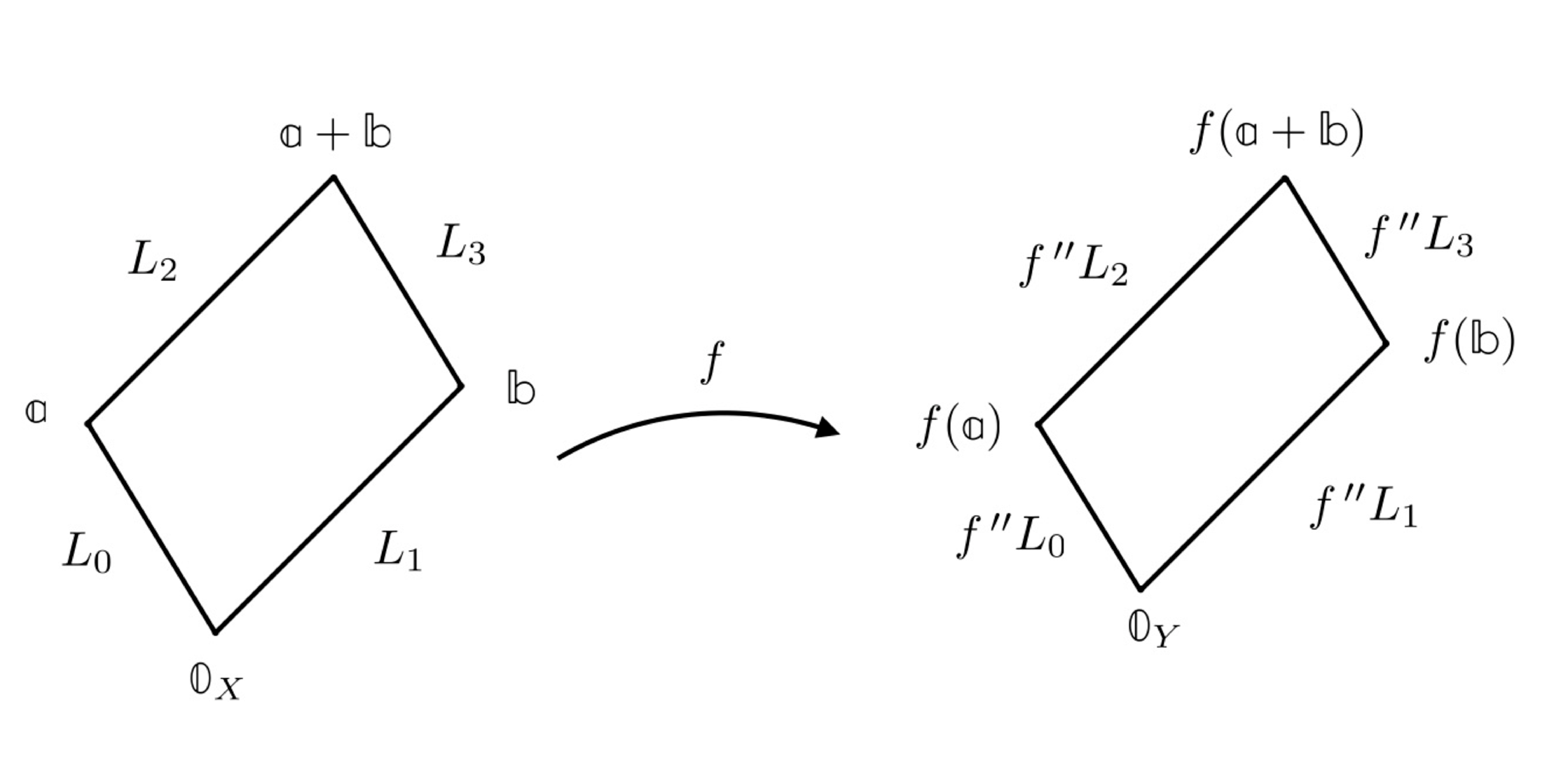}
\hfill\mbox{}

{\bf Case 2.} $\veca$ and $\vecb$ are not linearly independent. If one of $\veca$, $\vecb$  
equals to $\zerov_X$ then 
the additivity 
\begin{xitemize}
\xitem[lin-11] 
  $f(\veca+\vecb)=f(\veca)+f(\vecb)$
\end{xitemize}
trivially holds by \xitemof{lin-4}. 
Thus we may assume \wolog\ that $\veca\not=\zerov_X$, $\vecb\not=\zerov_X$, and
$\vecb\in\reals\veca$. 
If $f(\veca)=\zerov_Y$ or $f(\vecb)=\zerov_Y$, then $f\imageof\reals\veca=\ssetof{\zerov_Y}$ 
by \xitemof{lin-3}, \xitemof{lin-3-0} and \xitemof{lin-4} and the equation \xitemof{lin-11} 
again holds. Thus we may assume that $f(\veca)\not=\zerov_Y$, $f(\vecb)\not=\zerov_Y$ and
$f(\vecb)\in\reals f(\veca)$.

Let $\veca_0$, $\veca_1\in X$ be 
as in \xitemof{lin-8}. Then there is $i\in 2$ \st\ $f(\veca_i)$ and $f(\veca)$ are independent.
Let $E$ be the plane with $\veca_i$, $\veca$, $\zerov_X\in E$. By 
\Lemmaof{P-lin-5},\,\assertof{1}, $f\restr E$ is 1-1 and $f\imageof E$ is a plane with
$f(\veca_i)$, $f(\veca)$, $f(\vecb)$, $\zerov_Y\in f\imageof E$. 

Let $L_0=\veca\droite\zerov_X$, $L_1=(\veca_i+\veca)\droite\veca_i$,
$L_2=\veca_i\droite\zerov_X$, 
$L_3=(\veca_i+\veca)\droite\zerov_X$, $L_4=(\veca_i+\veca)\droite\veca$,
$L_5=(\veca_i+\veca+\vecb)\droite\vecb$, and $L_6=(\veca_i+\veca+\vecb)\droite(\veca+\vecb)$. 

We have $L_0\cap L_2\cap L_3=\ssetof{\zerov_X}$, $L_0\cap L_4=\ssetof{\veca}$,
$L_0\cap L_5=\ssetof{\vecb}$, $L_0\cap L_6=\ssetof{\veca+\vecb}$, 
$L_1\cap L_2=\ssetof{\veca_i}$, $L_1\cap L_3\cap L_4=\ssetof{\veca_i+\veca}$, and 
$L_1\cap L_5\cap L_6=\ssetof{\veca_i+\veca+\vecb}$. 

We also have that $L_0$ and $L_1$ are parallel to each 
other, $L_2$, $L_4$ and $L_6$ are parallel to each other, as well as $L_3$ and $L_5$ are 
parallel to each other. Since $f$ transfers this constellation keeping all the parallelisms, it follows, similarly to 
the Case 1., we can conclude in turn that $f(\veca_i+\veca)=f(\veca_i)+f(\veca)$,
$f(\veca_i+\veca+\vecb)=f(\veca_i)+f(\veca)+f(\vecb)$, and finally 
$f(\veca+\vecb)=f(\veca)+f(\vecb)$.
\bigskip

\noindent
\mbox{}\hfill
\includegraphics[width=0.8\textwidth]{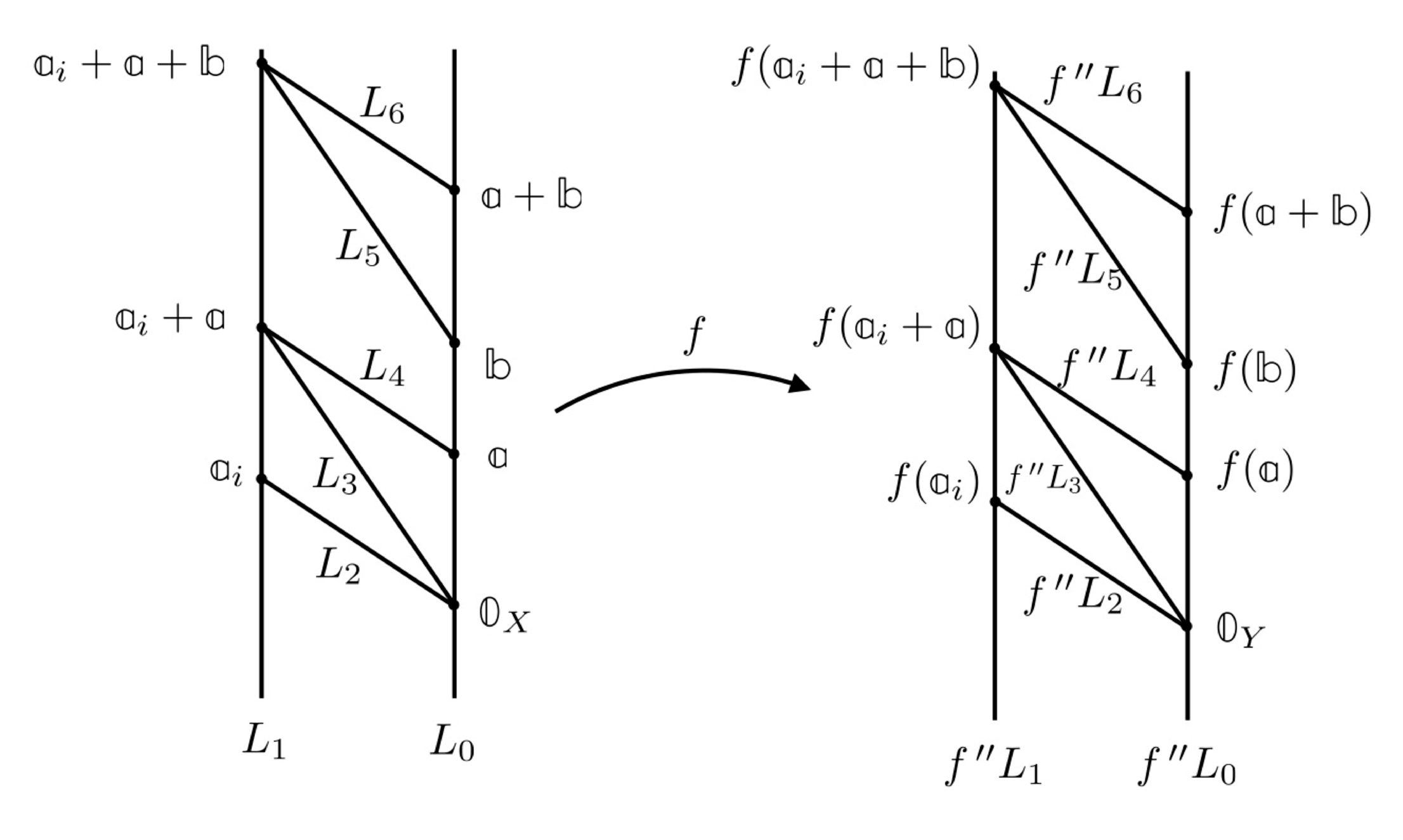}
\hfill\mbox{}

{\bf Case 3.} $\veca$ and $\vecb$ are independent in $X$ but $f(\veca)$ and $f(\vecb)$ are 
not independent in $Y$. If $f(\veca)=f(\vecb)=\zerov_Y$ then
$f(\veca+\vecb)=\zerov_Y=f(\veca)+f(\vecb)$ by the the Case 1 in the proof of 
\Lemmaof{P-lin-5},\,\assertof{1}. 
If one of $f(\veca)$ and $f(\vecb)$ is $\zerov_Y$, then 
we have $f(\veca+\vecb)=f(\veca)+f(\vecb)$ by \Lemmaof{P-lin-5-0}. Thus we may assume that  
$f(\veca)\not=\zerov_Y$ and $f(\vecb)\not=\zerov_Y$. Let $E\subseteq X$ be the plane with
$\zerov_X$, $\veca$, $\vecb\in E$. By the assumption of the present case, we have
$f\imageof E=\reals f(\veca)=\reals f(\vecb)$. Thus, there is $i\in 2$ \st\ $f(\veca_i)$ is 
independent from $f(\veca)$ (and also from $f(\vecb)$\,). 

In this case, some of the lines connecting $\zerov_X$ $\veca$, $\vecb$, $\veca_i$,
$\veca_i+\veca$, $\veca_i+\veca+\vecb$, $\veca+\vecb$ are sent to the same lines by $f$. 
However, parallelism of some lines survive to conclude in turn that $f(\veca_i+\veca)=f(\veca_i)+f(\veca)$,
$f(\veca_i+\veca+\vecb)=f(\veca_i)+f(\veca)+f(\vecb)$, and finally 
$f(\veca+\vecb)=f(\veca)+f(\vecb)$ just as in Case 2.
\bigskip

\noindent
\mbox{}\hfill
\includegraphics[width=0.8\textwidth]{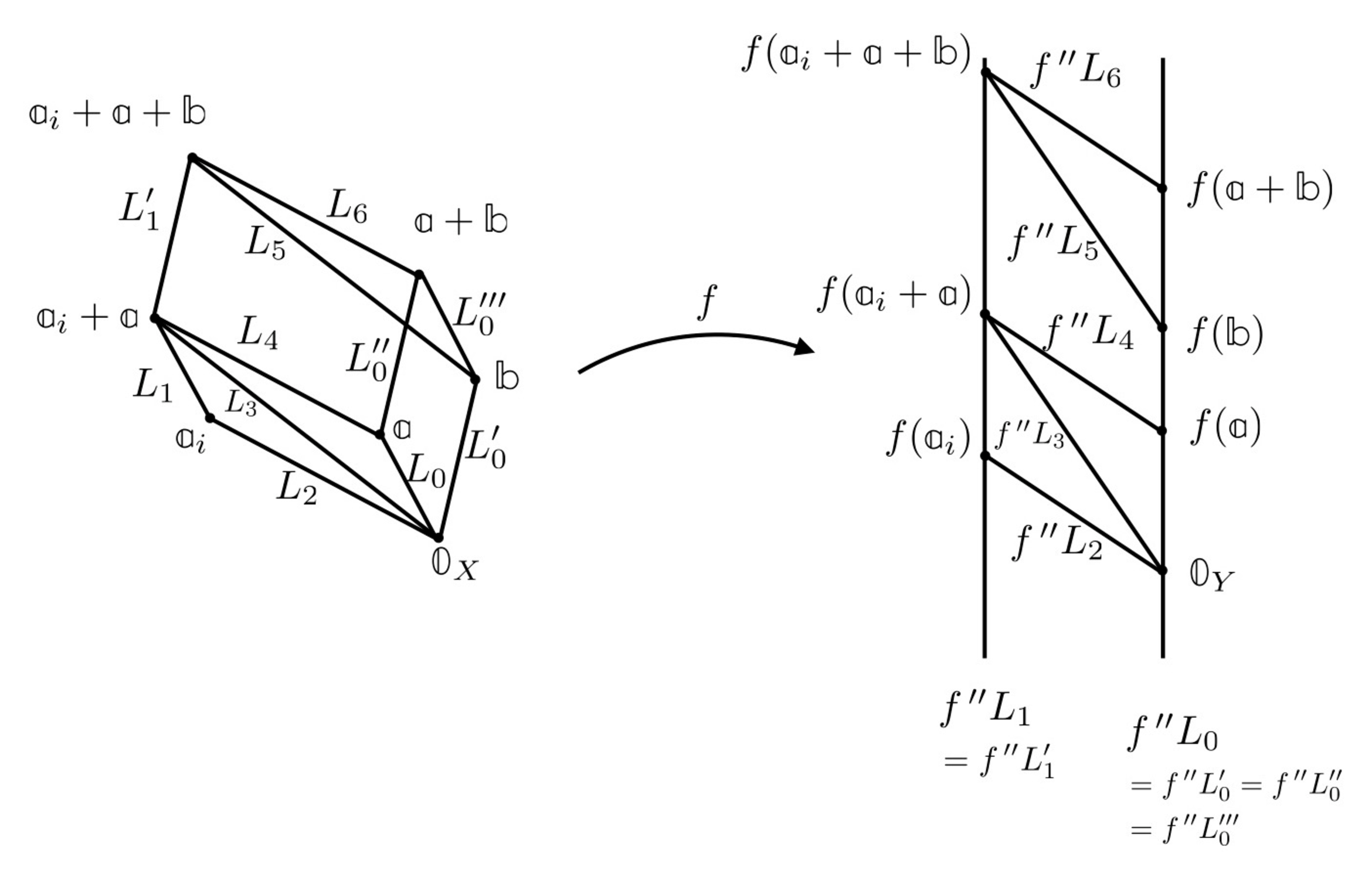}
\hfill\mbox{}

\qedofThm
\qedskip

\Thmabove\ can be easily generalized to a characterization of affine mappings
(i.e mappings $\mapping{g}{X}{Y}$ for $\reals$-linear spaces which can be represented as
$\veca\mapsto f(\veca)+\vecb$ for a linear mapping $\mapping{f}{X}{Y}$ and $\vecb\in Y$) 
with the following condition corresponding to \xitemof{lin-8}:

\begin{Cor}
  \Label{P-lin-7} Suppose that $X$, $Y$ are $\reals$-linear spaces and 
  $\mapping{g}{X}{Y}$ satisfies
  \begin{xitemize}
  \item[\xitemdof{lin-8}{'}\ \ \ ] there are $\veca^*$ $\veca^*_0$, $\veca^*_1\in X$ \st\
    $g(\veca^*_0)-g(\veca^*)$ and $g(\veca^*_1)-g(\veca^*)$ are linearly independent.      
  \end{xitemize}

  Then $g$ is an affine mapping if and only if $g$ satisfies
  \begin{xitemize}
  \item[\xitemof{lin-3}\ \ \ ] the image $g\imageof L$ of any line $L\subseteq X$ by $g$ 
    is either a point or a line in $Y$; and 
  \item[\xitemof{lin-7}\ \ \ ] for any line $L\subseteq X$, if $g\imageof L$ is also a line in $Y$, then
    $g\restr L$ is 1-1.
  \end{xitemize}
\end{Cor}
\prf If $\mapping{g}{X}{Y}$ is an affine mapping then $g$ clearly satisfies
\xitemof{lin-3} and \xitemof{lin-7}.  

Conversely, suppose that $\mapping{g}{X}{Y}$ satisfies \xitemdof{lin-8}{'}, \xitemof{lin-3} 
and \xitemof{lin-7}. 
Let $\veca^*$, $\veca^*_0$, $\veca^*_1\in X$ be as in \xitemdof{lin-8}{'}.
Let $\mapping{f}{X}{Y}$ be defined by 
\begin{xitemize}
\xitem[lin-11-0] 
  $f(\veca)=g(\veca+\veca^*)-g(\veca^*)$ for
  $\veca\in X$.
\end{xitemize}

Clearly $f$ also satisfies \xitemof{lin-3} and $\xitemof{lin-7}$. We have 
$f(\zerov_X)=g(\veca^*)-g(\veca^*)=\zerov_Y$.
$f(\veca^*_i-\veca^*)=g(\veca^*_i-\veca^*+\veca^*)-g(\veca^*)=g(\veca^*_i)-g(\veca^*)$ for
$i\in 2$. Thus, letting $\veca_i=\veca^*_i-\veca^*$ for $i\in 2$, we see that $f$ satisfies 
\xitemof{lin-8}.
By \Thmof{P-lin-6}, it follows that $f$ is a linear mapping. Since we have 
\begin{xitemize}
\xitem[] 
  $g(\veca)=g((\veca-\veca^*)+\veca^*)=f(\veca-\veca^*)+g(\veca^*)
  =f(\veca)+f(-\veca^*)+g(\veca^*)=f(\veca)+g(\zerov_X)$ 
\end{xitemize}
for all $\veca\in X$,  
$g$ is an affine mapping.\qedofCor\qedskip

The following \xitemof{lin-12} apparently holds for any affine mapping $\mapping{f}{X}{Y}$.
Conversely, it is easy to see that \xitemof{lin-3} and \xitemof{lin-12} imply
\xitemof{lin-3-0}. Thus we obtain: 

\begin{Cor}
  \Label{P-lin-8}
  Suppose that $X$, $Y$ are $\reals$-linear spaces and 
  $\mapping{g}{X}{Y}$ satisfies
  \begin{xitemize}
  \item[\xitemdof{lin-8}{'}\ \ \ ] there are $\veca^*$ $\veca^*_0$, $\veca^*_1\in X$ \st\
    $g(\veca^*_0)-g(\veca^*)$ and $g(\veca^*_1)-g(\veca^*)$ are linearly independent.      
  \end{xitemize}
  Then $g$ is an affine mapping if and only if $g$ satisfies 
  \begin{xitemize}
  \item[\xitemof{lin-3}\ \ \ ] the image $g\imageof L$ of any line $L\subseteq X$ by $g$ 
    is either a point or a line in $Y$; and,  
  \xitem[lin-12] for any $\veca$, $\vecb$, $\vecc\in X$ with
    $\vecc\in\interval\concat\vecb$, $g(\veca)$, either $g(\veca)=g(\vecb)=g(\vecc)$, or 
    $g(\vecc)\in g(\veca)\concat g(\vecb)$ holds.  \qed
  \end{xitemize}
\end{Cor}

The following characterization of linear mappings consisting of a mixture of algebraic and 
geometric properties is also easy to prove\ifextended\else\ as an application of \Lemmaof{P-lin-4}\fi.

\begin{Prop}
  \Label{P-lin-9} Suppose that $X$ and $Y$ are $\reals$-linear spaces 
  and $\mapping{f}{X}{Y}$. Then $f$ is linear if and only if
  \begin{xitemize}
  \xitem[] $f$ is an additive function; and 
  \item[\xitemdof{lin-12}{'}\ \ \ ] for any $\veca$, $\vecc\in X$, if
    $\vecc\in\veca\concat\zerov_X$, then, either $f(\veca)=f(\vecc)=\zerov_Y$ or 
    $\vecc\in f(\veca)\concat\zerov_Y$. \ifextended\else\qed\fi
  \end{xitemize}
\end{Prop}
\ifextended{\tt
\prf For $\veca\in X$, let $\mapping{\varphi_\veca}{\reals}{\reals}$ be \st\
\begin{xitemize}
\xitemA[lin-13] 
  $f(r\veca)=\varphi_\veca(r)f(\veca)$ for all $r\in\reals$. 
\end{xitemize}
We assume, \wolog, that
$\varphi_\veca(r)=r$ for all $r\in\reals$, if $f(\veca)=\zerov_Y$. 
\begin{Claim}
  $\varphi_\veca$ is additive and monotonically increasing.
\end{Claim}
\prfofClaim 
If $f(\veca)=\zerov_Y$, the assertion is trivial. So assume that $f(\veca)\not=\zerov_Y$. 

For $r$, $s\in\reals$, we have
$\varphi_\veca(r+s)f(\veca)=f((r+s)\veca)=f(r\veca+s\veca)=f(r\veca)+f(s\veca)
=\varphi_\veca(r)f(\veca)+\varphi_\veca(s)f(\veca)=(\varphi_\veca(r)+\varphi_\veca(s))f(\veca)$.
This implies that $\varphi_\veca(r+s)=\varphi_\veca(r)+\varphi_\veca(s)$.

To show that $\varphi_\veca$ is monotonically increasing,
suppose that $r$, $s\in\reals$ with $r<s$. For simplicity, consider the case $1<r<s$. All 
other cases can be treated similarly. Since $\veca\in r\veca\concat\zerov_X$, we have
$f(\veca)\in f(r\veca)\concat\zerov_Y$ by \xitemdof{lin-12}{'}. By definition 
\xitemAof{lin-13}, it follows that $1<\varphi_\veca(r)$. Since
$r\veca\in s\veca\concat\zerov_X$, it follows that $\varphi_\veca(r)<\varphi_\veca(s)$. 
\qedofClaim\qedskip

By \Lemmaof{P-lin-4}, it follows that $\varphi_\veca$ is $\reals$-linear since
$\varphi_\veca(1)=1$, $\varphi_\veca$ must be the identity function. \qedofProp
}\fi 

\end{document}